\documentclass[12pt]{article}

\usepackage[latin1]{inputenc}
\usepackage[T1]{fontenc}
\usepackage[english]{babel}
\usepackage{lmodern}
\usepackage[babel=true,kerning=true]{microtype}
\usepackage{amsmath,amssymb}
\usepackage{amsfonts}
\usepackage{amsthm}
\usepackage{appendix}
\usepackage{mathrsfs}
\usepackage{tabularx}
\usepackage{pgf,tikz}
\usetikzlibrary{arrows}
\usepackage{vmargin}
\usepackage{variations}
\usepackage{tkz-tab}
\usepackage[normalem]{ulem}
\usepackage{ulem}
\usepackage{eurosym}
\usepackage{fancyhdr}
\usepackage{latexsym}
\usepackage{bbding}
\usepackage{vwcol}
\usepackage{multirow}
\usetikzlibrary{decorations.pathreplacing}

\theoremstyle{plain}

\theoremstyle{definition}

\theoremstyle{remark}

\newcommand{\I}{\mathop{}\mathopen{}\ensuremath{\mathrm{i}}\mathop{}\mathopen{}}

\setmarginsrb{2cm}{2cm}{2cm}{2cm}{0cm}{0cm}{0cm}{0cm}

\title{A direct proof of the irrationality of $\tan^2(r \pi)$}
\author{Lionel Ponton}
\date{}
\begin{document}

\maketitle
\thispagestyle{empty}

\begin{abstract} Given a rational number $r$ such that $2r$ is not an integer, we prove that $\tan^2(r\pi)$ is irrational unless it is equal to $0$, $1$, $3$ or $\frac{1}{3}$, using only basic trigonometry and the Rational Root Theorem. Moreover, we deduce that $\tan(r\pi)$, $\cos^2(r\pi)$ ans $\cos(r\pi)$ are irrational numbers except in usual cases.
\end{abstract}

Let $r=\frac{d}{n}$ be a rational number written in lowest terms such that $n$ is not equal to $2$. It is well-known that $\tan(r \pi)$ is an irrational number unless it is equal to $-1$, $0$ or $1$ and this result is usually deduced from a similar property of $\cos(r \pi)$ (see, e.g., \cite[Th. 6.16]{NZM91}).  The aim of this short note is to give a direct proof using only basic trigonometry and the Rational Root Theorem \cite[Th. 6.13]{NZM91} which we recall here:

\medskip

\textit{Rational Root Theorem}. If a polynomial with integral coefficients $P=\sum_{k=0}^{n} c_kX^k$ has a nonzero rational root $\frac{a}{b}$ where the integers $a$ and $b$ are relatively prime, then $a$ divides $c_0$ and $b$ divides $c_n$. 

\medskip

Note that in \cite[Th. 5, p. 95]{ES03}, the authors give a direct proof in the case $d=1$ due to P. Tur\'an but, even if it uses the same tools, the following is shorter, requires fewer calculations and provides a more general statement. Indeed, we not only prove the irrationality of $\tan(r \pi)$ in suitable cases but we show that $\tan^2(r \pi)$ is an irrational number unless it is equal to $0$, $1$, $3$ or $\frac{1}{3}$. Moreover, this allows us to derive that $\cos^2(r\pi)$ and $\cos(r\pi)$ are irrational numbers except in usual cases.

Assume that $\tan^2(r\pi)$ is a rational number not equal to $0$ or $1$.  Note that 
\begin{equation}
\label{eq_double}
\tan^2(2r \pi) = \frac{4\tan^2(r \pi)}{(1-\tan^2(r \pi))^2}
\end{equation}
	is also a rational number not equal to $0$ or $1$. Indeed, clearly $\tan^2(2r \pi)$ is a nonzero rational number and $\tan^2(2r \pi) \neq 1$ since the solutions of $\frac{4x}{(1-x)^2}=1$ are the irrational numbers $x=3\pm 2\sqrt{2}$. Thus, possibly multiplying by a suitable power of 2, we can assume that $n$ is odd and write $n=2m+1$ where $m$ is a positive integer (since $\tan^2(r \pi) \neq 0$).

Let us put $t=\tan(2r \pi)=\tan(\frac{2d}{n}\pi)$ and $s=t^2=\tan^2(2r \pi)$. Then, we have
\[(1+\I t)^n=\left(\frac{\mathrm{e}^{\I\frac{2d}{n}\pi}}{\cos(\frac{2d}{n}\pi)}\right)^n=\frac{1}{\cos^n(\frac{2d}{n}\pi)}=\left(\frac{\mathrm{e}^{-\I\frac{2d}{n}\pi}}{\cos(-\frac{2d}{n}\pi)}\right)^n=(1-\I t)^n,\]
so $t$ is a zero of
\[P:=(-1)^m\frac{(1+\I X)^n-(1-\I X)^n}{2\I X}=\sum_{j=0}^m \binom{n}{2j+1} (-1)^{m+j} X^{2j}.\]
Hence, $s$ is a nonzero root of the monic polynomial with integral coefficients
\[Q:=\sum_{j=0}^m (-1)^{m+j}\binom{n}{2j+1} X^{j}.\]
Since $s$ is a rational number, we conclude by the Rational Root Theorem that $s$ is an integer and $s$ divides $(-1)^mn$ so $s$ is odd.

Furthermore, since $s$ is rational, the previous argument applied to $u:=\tan^2(4r \pi)$ instead of $s$ ensures that $u$ is an odd integer and, by applying \eqref{eq_double} to $2r$ instead of $r$, we see that $s$ is a solution of the quadratic equation $(E_u):u x^2-2(u+2)x+u=0$. Thus, since $s\geqslant 1$,
\begin{equation}
\label{eq_expression_of_s}
s=\frac{u+2+2\sqrt{u+1}}{u},
\end{equation}
and so $k:=\sqrt{u+1}$ is rational. Moreover, since $k$ is a nonzero root of $X^2-(u+1)$, by applying the Rational Root Theorem once again, we deduce that $k$ is a positive integer. Thereby, due to equation \eqref{eq_expression_of_s},
\[s=\frac{k^2+1+2k}{k^2-1}=\frac{(k+1)^2}{k^2-1}=\frac{k+1}{k-1}.\]
Furthermore, since $u=(k-1)(k+1)$ is odd, $k-1$ and $k+1$ are coprime. Since $s$ is an integer, it follows that $k=2$ and thus $s=3$. But, similarly, by \eqref{eq_double}, $\tan^2(r\pi)$ is a solution of $(E_3):3x^2-10x+3=0$, and so $\tan^2(r\pi)=\frac{1}{3}$ or $\tan^2(r\pi)=3$.

We conclude that $\tan^2(r\pi)$ is irrational unless it is equal to $0=\tan^2(0)$, $1=\tan^2(\frac{\pi}{4})$, $\frac{1}{3}=\tan^2(\frac{\pi}{6})$ or $3=\tan^2(\frac{\pi}{3})$. 

Lastly, since $\sqrt{3}$ is irrational, it follows that $\tan(r \pi)$ is an irrational number unless it is equal to $-1$, $0$ or $1$. Moreover, since $1+\tan^2(r \pi)=\frac{1}{\cos^2(r \pi)}$, we also conclude that $\cos^2(r \pi)$ is an irrational number unless it is equal to $0$, $1$, $\frac{1}{2}$, $\frac{1}{4}$, $\frac{3}{4}$ and thus $\cos(r\pi)$ is irrational unless it is equal to $0$, $\pm 1$ or $\pm \frac{1}{2}$.

\bigskip

\textbf{Acknowledgement.} The author wishes to express his gratitude to Daniel Duverney for pointing out that the earlier version of the proof, which only treated the irrationality of $\tan(r \pi)$, could easily be adapted to $\tan^2(r\pi)$.

\end{document}